\documentclass[11pt]{amsart}

\input xy
\xyoption{all}

\begin{document}

\date{November 19, 2003}

\title[Homology of Linear Groups]{Homology of Linear Groups via Cycles in
$BG\times X$}

\author{Kevin P.~Knudson}\thanks{First author partially supported by NSF grant nos. 
DMS-0070119, DMS-0242906.}
\address{Department of Mathematics and Statistics, Mississippi State University,
P.~O.~Drawer MA, Mississippi State, MS 39762}
\email{knudson@math.msstate.edu}

\author{Mark E.~Walker}\thanks{Second author partially supported by the NSA}
\address{Department of Mathematics, University of 
Nebraska--Lincoln, Lincoln, Nebraska 68588}
\email{mwalker@math.unl.edu}

\newtheorem{theorem}{Theorem}[section]
\newtheorem{prop}[theorem]{Proposition}
\newtheorem{lemma}[theorem]{Lemma}
\newtheorem{cor}[theorem]{Corollary}
\newtheorem{conj}[theorem]{Conjecture}
\newtheorem{definition}[theorem]{Definition}
\newtheorem{remark}[theorem]{Remark}
\newtheorem{guess}[theorem]{Equivariant Isomorphism Conjecture}

\newcommand{\zz}{{\mathbb Z}}
\newcommand{\zq}{{\mathbb Q}}
\newcommand{\gm}{{\mathbb G}_m}
\newcommand{\ga}{{\mathbb G}_a}
\newcommand{\ra}{\rightarrow}
\newcommand{\lra}{\longrightarrow}
\newcommand{\bop}{\bigoplus}
\newcommand{\hig}{{\mathcal H}_i}
\newcommand{\hih}{{\mathcal H}_i}
\newcommand{\hgi}{{\mathcal H}^i}
\newcommand{\zl}{{\mathbb Z}/\ell}
\newcommand{\zn}{{\mathbb Z}/n}
\newcommand{\gln}{GL_n}
\newcommand{\glnk}{GL_n(k)}
\newcommand{\rbar}{\overline{R}}
\newcommand{\aff}{{\mathbb A}^1}
\newcommand{\kbar}{\overline{k}}
\newcommand{\glm}{GL_m}
\newcommand{\glmk}{GL_m(k)}
\newcommand{\zr}{{\mathbb R}}
\newcommand{\zc}{{\mathbb C}}
\newcommand{\R}{{\mathbb R}}
\newcommand{\C}{{\mathbb C}}
\newcommand{\bP}{{\mathbb P}}

\newcommand{\abar}{\underbar{A}}
\newcommand{\sbar}{\overline{S}}
\newcommand{\cH}{{\mathcal H}}
\newcommand{\cM}{{\mathcal M}}
\newcommand{\chr}{\text{char}}
\newcommand{\cF}{\mathcal F}
\newcommand{\cA}{{\mathcal A}}

\newcommand{\cO}{{\mathcal O}}
\newcommand{\bark}{{\overline{k}}}
\newcommand{\Rbar}{{\overline{R}}}
\newcommand{\Deldot}{\Delta^\bullet}
\newcommand{\Tot}{\operatorname{Tot}}
\newcommand{\Hom}{\operatorname{Hom}}
\newcommand{\Spec}{\operatorname{Spec}}

\maketitle

\section{Introduction}

Homology theories for algebraic varieties are often constructed using
simplicial sets of algebraic cycles. For example, Bloch's higher Chow
groups and motivic cohomology, as defined by Suslin and Voevodsky
\cite{VSF}, 
are given in this fashion.  In this paper, we construct homology
groups $\hig(X,G)$, where $G$ is an algebraic group and $X$ is a
variety,
by considering cycles on the simplicial scheme $BG \times X$,
an idea first suggested by Andrei Suslin.
If $X= \Spec(R)$ is an affine scheme, then there is a natural map
$$
H_i(G(R),A)\lra \hig(X,G;A),
$$ 
whose source is the usual group homology of the discrete group $G(R)$
of $R$-points of the algebraic group $G$. Moreover, this map is an
isomorphism if $R = k$ and $k$ is algebraically closed,
so that these groups capture 
the homology of the discrete group $G(k)$.  The functors $\hig(-,G;A)$
are naturally equipped with transfer maps in the sense of
Suslin--Voevodsky \cite{susvoe}.  We also have cohomological versions
$\hgi(X,G;A)$ and by applying these to the standard cosimplicial
variety $\Delta^\bullet_k$ we deduce a spectral sequence
$$
E_1^{s,t} = {\mathcal H}^s(\Delta^t,G;\zn) \Longrightarrow
H_{\text{\'et}}^{s+t}(BG_k,\zn)
$$
for $k$ algebraically closed with $\text{char}(k)$ not dividing $n$.
Since $\hgi(\Delta^0,G; A)
\cong H^i(G(k), A)$, the bottom row of the spectral sequence
is just $H^\bullet(G(k),\zn)$.  Thus, this spectral sequence
provides a mechanism for comparing the group cohomology of the
discrete group $G(k)$ with 
the \'etale cohomology of the simplicial scheme $BG_k$.

We recall the following conjecture of E.~Friedlander \cite{friedmis}.

\begin{conj}[Friedlander's Generalized Isomorphism Conjecture]
If $G$ is an algebraic group over an
algebraically closed field $k$, then the natural map of simplicial
schemes $BG(k)\ra BG_k$ (with the first one being ``discrete'') 
induces an isomorphism
$$
H_{\text{\em \'et}}^\bullet(BG_k,\zn)\lra 
H^\bullet(BG(k),\zn) = H^\bullet(G(k), \zn)
$$ 
provided $\text{\em char}(k)$ does not divide $n$.  
\end{conj}

If $k={\mathbb C}$, then the Isomorphism Conjecture
reduces to the claim that the identity map $G^\delta\ra
G^{\text{top}}$, where $G^\delta$ is the Lie group $G^{\text{top}}$
viewed as a discrete group, induces an isomorphism in cohomology with
finite coefficients:
$$
H^\bullet(BG^{\text{top}},\zn)
\ra 
H^\bullet(BG^\delta,\zn) =
H^\bullet(G^\delta,\zn).
$$  
In this
context, it also makes sense to consider the field $k={\mathbb R}$ and
assert this isomorphism for real Lie groups.

The Isomorphism Conjecture is known to hold in the following cases:

\begin{enumerate}
\item $k={\overline{\mathbb F}}_p$, $G$ arbitrary \cite{friedmis};

\item $G$ solvable \cite{friedmis},\cite{jardine},\cite{milnor};

\item $G$ is one of 
the stable groups $GL,SL,Sp,O$ \cite{jardine},\cite{karoubi},\cite{suslinlocal};

\item $H^2(G,\zn)$ for $G$ a real Chevalley group \cite{sahwag};

\item $H^3(GL_m(k),\zn)$ ($k$ arbitrary), $H^3(SL_2({\mathbb
  C}),\zn)$, \cite{knudson},\cite{sah},\cite{suslinbloch}. 
\end{enumerate}

Using the spectral sequence described above, we are able to deduce all
the known cases of the Isomorphism Conjecture immediately (except for
the results about {\em real} Lie groups).  Moreover, we can explicitly
identify the group $H_{\text{\'et}}^4(BGL_{m},\zn)$.

\medskip

\noindent {\bf Corollary \ref{h4}.} {\em There is an exact sequence
$$0
\lra H_{\text{\em \'et}}^4(BGL_{m},\zn)\lra H^4(B\glmk,\zn)\stackrel{d}{\lra} 
{\mathcal H}^4(\Delta^1,GL_m;\zn)
$$ where $d$ is the difference of the maps
induced by the two face maps $\Delta^0 \ra \Delta^1$.}

\medskip

\noindent Conjecturally, the map $d$ is the zero map, of course.

This paper is organized as follows.  In Section \ref{functors} we
construct the functors $\hig(X,G;A)$ and discuss their basic
properties.  Section \ref{examples} contains some relatively easy
examples.  The spectral sequence is constructed in Section \ref{ss}.
There is a relationship between the $\hig(\text{Spec}(k),G;A)$ and
the homology groups $H_i(BG(\kbar)/\Gamma;A)$,
where $\Gamma=\text{Gal}(\kbar/k)$. We discuss this in Section
\ref{invariant} and prove the following results.

\medskip

\noindent {\bf Theorem \ref{specr}.}  {\em Let $G$ be an algebraic
group over $\zr$ and let $A=\zq$ or $\zn$ where $n$ is odd.  Then
there are canonical isomorphisms
$$
\hig(\Spec(\zr),G;A)\cong H_i(G(\zc),A)^{{\mathbb Z}/2} 
$$ 
for all $i\ge 0$.}

\medskip

If $A=\zz/2$, then the structure of $\hig(\Spec(\zr),G;A)$ is
much more complicated.  However, we can deal with the case $G=\gm$.

\medskip

\noindent {\bf Theorem \ref{specrgm}.} {\em For all $k\ge 1$,
$$\cH_{2k}(\Spec(\zr),\gm;\zz/2) \cong
\cH_{2k+1}(\Spec(\zr),\gm;\zz/2) \cong
(\zz/2)^{k-1}.$$}

This is proved by noting that an equivariant version of Friedlander's
conjecture holds for $\gm$, and then using the calculation of the
Bredon cohomology of $B\gm(\zc)^{\text{top}}\simeq \zc \bP^\infty$ given
in \cite{dugger}.

We are also able to handle unipotent groups over $\zr$.  We discuss this
at the end of Section \ref{invariant}.

\medskip

\noindent {\em Notation.} Throughout, $k$ denotes a field and $p$
denotes the exponential characteristic of $k$; that is, $p=1$ if
$\text{char}(k) = 0$ and $p=\text{char}(k)$ if $\text{char}(k)>0$.

\medskip

\noindent {\em Acknowledgments.} The construction of the
$\hig(X,G;A)$ was suggested to us by Andrei Suslin.  We thank him for
putting us on this path.

\medskip

\noindent {\em Conventions:}
Throughout this paper, given a field $k$, a $k$-scheme is a 
separated scheme of finite type over $k$.

\section{Construction of the functors $\hig(X,G;A)$}\label{functors}
Let $G$ be a linear algebraic group over $k$.  Consider the simplicial
classifying scheme $BG$:
$$
\xymatrix{
\ast & G \ar@<-1ex>[l]\ar@<1ex>[l] & G^2 \ar@<-1ex>[l]\ar[l]
\ar@<1ex>[l] & G^3 \ar@<-1.5ex>[l]\ar@<1.5ex>[l]_\vdots
  \cdots & G^n \ar@<-1.5ex>[l]\ar@<1.5ex>[l]_(.4)\vdots & 
\cdots \ar@<-1.5ex>[l]\ar@<1.5ex>[l]_(.4)\vdots}.
$$
For any $k$-scheme $X$, 
we may form the product $BG\times X$:
$$\xymatrix{
X & G\times X \ar@<-1ex>[l]\ar@<1ex>[l] & G^2\times X\ar@<-1ex>[l]\ar[l]
\ar@<1ex>[l] & G^3 \times X\ar@<-1.5ex>[l]\ar@<1.5ex>[l]_\vdots
  \cdots & G^n\times X\ar@<-1.5ex>[l]\ar@<1.5ex>[l]_(.4)\vdots & 
\cdots\ar@<-1.5ex>[l]\ar@<1.5ex>[l]_(.4)\vdots}.
$$
If $U$ and $V$ are smooth $k$-schemes, define $C_0(U\times V/V)$ as
$$
C_0(U\times V/V)=\zz\left\{\begin{array}{l}
{\text{closed integral subschemes}\,Z\subset U\times V\,\text{such}} \\
\text{that}\, Z\lra V\,\text{is finite and surjective over} \\
\text{some connected component of}\, V \end{array}\right\}.
$$
This construction was first studied by Suslin-Voevodsky in
\cite{susvoe}.

The functor $C_0(U\times V/V)$ is covariant in $U$ via pushforward of
cycles. In more detail,
if $p:U\ra U'$ is a morphism of $k$-schemes and if $Z\subset U\times V$ is
a closed integral subscheme that is finite and surjective over some
connected component of $V$, then $Z' = (p\times \text{id}_V)(Z)$ is
finite and surjective over some component of $V$.
The map $p_*$ is defined by sending $[Z]$ (the cycle
associated to $Z$) to $d [Z']$, where $d$ is the
degree of the finite, dominant map $Z \to Z'$, and then extending
$p_*$ to all of $C_0(U \times V/V)$ by linearity.

The functor $C_0(U \times V/V)$ is contravariant in $V$ via pullback
of cycles. Given a morphism $f: V' \to V$ of smooth $k$-schemes, for
a closed, integral
subscheme $Z$ of $U \times V$ that is finite and surjective over a
component of $V$, the pullback $f^{-1}(Z) = V' \times_V Z$ is a closed
subscheme of $U \times V'$ each integral component of which is
finite and surjective over a component of $V'$. One defines
$f^*([Z])$ to be the cycle in $U \times V'$ associated to $f^{-1}(Z)$
by taking multiplicities of its integral components in the
usual fashion. The definition of $f^*$ is then extended to all of
$C_0(U \times V/V)$ by linearity.

In particular, we apply the covariant functor $C_0(- \times X/X)$
degreewise to the simplicial scheme $BG$ to obtain a simplicial
abelian group $C_0(BG\times X/X)$:
$$
\xymatrix{
C_0(\ast\times X/X) & C_0(G\times X/X) \ar@<-1ex>[l]\ar@<1ex>[l] 
 &\cdots\ar@<-1.5ex>[l]\ar@<1.5ex>[l]_(.3)\vdots
  & C_0(G^n\times X/X)\ar@<-1.5ex>[l]\ar@<1.5ex>[l]_(.7)\vdots & 
\cdots\ar@<-1.5ex>[l]\ar@<1.5ex>[l]_(.3)\vdots}.
$$

\begin{definition} Let $A$ be an abelian group. The group
  $\hig(X,G;A)$ is defined to be the $i$-th homotopy 
group of the simplicial abelian group $C_0(BG\times X/X)\otimes A$; that is,
$$
\hig(X,G;A)=h_i(C_0(BG\times X/X)\otimes A,d),
$$ 
where $d=\sum (-1)^kp_{k*}$ and
$p_k:G^i\ra G^{i-1}$ is the map 
$$
p_k(g_0,g_1,\dots ,g_{i-1}) = \begin{cases}
                                   (g_1,\dots ,g_{i-1}) & k=0 \\ 
                                   (g_0,\dots ,g_{k-1}g_k,\dots
  g_{i-1}) & 0<k<i-1 \\ 
                                   (g_0,\dots ,g_{i-2}) & k=i-1.
                                 \end{cases}
$$
We also define groups $\hgi(X,G;A)$ by
$$
\hgi(X,G;A) = h^i(\text{\em Hom}(C_0(BG\times X/X),A)).
$$
\end{definition}

\subsection{Variance} 
The notation $\hig(X,G; A)$ was chosen to suggest the same variance in $X$ and
$G$ as the bi-functor $\text{Hom}(-,-)$. Namely,
if $f:G\ra H$ is a morphism of linear algebraic groups over $k$, then
there is an induced map $f_*:\hig(X,G;A)\ra \hih(X,H;A)$ for any
scheme $X$.  Indeed, the push-forward maps $f^{\times n}_* : 
C_0(G^n \times X/X) \to C_0(H^n \times X/X)$, $n \geq 0$, are
compatible with the simplicial structures so that we obtain 
a map of simplicial abelian groups $f_*: C_0(BG \times X/X) \to C_0(BH
\times X/X)$ and hence a map on homotopy groups as desired.
Likewise, if $\varphi:X\ra Y$ is a morphism of smooth $k$-schemes, 
then pull-back of cycles determines maps
$\varphi^*: C_0(G^n \times X/X) \to C_0(G^n \times Y/Y)$, $n \geq 0$,
that are again compatible with simplicial structures so that we have
an induced map 
$\varphi^*:\hig(Y,G;A)\ra\hig(X,G;A)$ on homotopy groups.

Of course, we have the opposite variance for the functors
$\hgi(X,G;A)$ ---
these are contravariant in $G$ and covariant in $X$.

\subsection{Sheaf-theoretic interpretation of $C_0(BG\times X/X)$}
We now present an equivalent definition, due to Suslin-Voevodsky,
\cite{susvoe} 
of $\hig(X,G)$ in
terms of sheaves in the qfh topology.  Recall \cite{susvoe} that a
morphism $q:X\ra Y$ is a {\em topological epimorphism} if the underlying
Zariski topological space of $Y$ is a quotient space of the underlying
Zariski topological space of $X$.  The map $q$ is a {\em universal
topological epimorphism} if for any $Z\ra Y$ the morphism $q_Z:X\times_Y Z
\ra Z$ is a topological epimorphism.
An h-covering of a scheme $X$ is a finite family of morphism of finite
type $\{p_i:X_i\ra X\}$ such that $\coprod p_i:\coprod X_i\ra X$ is a
universal topological epimorphism.  A qfh-covering of $X$ is an h-covering
$\{p_i\}$ such that all the morphisms $p_i$ are quasi-finite.

If $X$ is a smooth separated scheme of finite type over $k$, let 
$\zz[1/p]_{\text{qfh}}(X)$ denote
the sheaf in the qfh topology associated to
the presheaf 
$$
T\mapsto \zz[1/p]\text{Hom}_{\text{Sch}/k}(T,X),
$$
where $\zz[1/p] S$ denotes the free $\zz[1/p]$-module on
the set $S$ and 
$p$ is the exponential characteristic of $k$.

\begin{theorem}[\cite{susvoe},6.7]\label{sheaf}  If $Y$ is a separated
scheme, then there is an isomorphism 
$$
C_0(Y \times X/X)\otimes \zz[1/p] \cong 
\Gamma(X,\zz[1/p]_{\text{\em qfh}}(Y))
$$
that is natural in both $X$ and $Y$.
\end{theorem}

Note that this theorem implicitly asserts that $C_0(Y\times -/-)[1/p]$
is a sheaf in the qfh-topology. As shown in \cite[\S 5]{susvoe}, any
sheaf ${\mathcal F}$ in the qfh topology admits transfer maps; that is, for
any finite surjective map $f:X \ra S$, where $X$ is reduced and irreducible
and $S$ is irreducible and regular, there is a transfer homomorphism
$$
\text{Tr}_{X/S}:{\mathcal F}(X)\lra {\mathcal F}(S)
$$
satisfying certain expected properties (cf. 
\cite[4.1]{susvoe}). In the case $\cF = C_0(Y \times -/-)$, the
transfer homomorphism is defined by pushforward of cycles in the
evident manner.

By Theorem \ref{sheaf}, if $p$
is invertible in $A$, the complex $(C_0(BG\times -/-)\otimes A,d)$ is
a complex of qfh sheaves and so is equipped with degreewise transfer
maps.  These commute with $d$ and hence for any finite surjective map
$X\ra S$, we have transfer maps
$$\text{Tr}_{X/S}:\hig(X,G;A)\lra \hig(S,G;A).$$  Thus we have proved the
following proposition.

\begin{prop}\label{transfers} Suppose that $p$ is invertible in the abelian
group $A$.  Then the functor $\hig(-,G;A)$ is a presheaf with transfers. \hfill $\qed$
\end{prop}

\subsection{First examples}\label{examples}

\begin{prop} \label{finiteprop}
Let $G$ be a finite group. For all $i\ge 0$, $\hig(X,G;\zz)=
H_i(G,\zz)$ for any integral scheme $X$.
\end{prop}

\begin{proof} An irreducible subscheme $Z\subset G^i\times X$ that is
finite and surjective over $X$ must be isomorphic to $X$ since the scheme
$G^i\times X$ is simply the disjoint union of copies of $X$.  It follows
that the group $C_0(G^i\times X/X)$ is the free abelian group on the
closed points of $G^i$.  Thus, $\hig(X,G;A)$ is the $i$-th homology of the
standard complex for computing $H_\bullet(G,\zz)$.
\end{proof}

In this case the pullback and transfer maps are easy to describe.  
For any morphism $f: X \ra S$, 
$f^*:\hig(S,G;\zz)\ra\hig(X,G;\zz)$ is induced by pullback of cycles,
which in this case is the map sending 
$[\{g\} \times S]$ to $[\{g\} \times X]$. Thus
  $f^*$ is simply the identity map (after making the identification of
  Proposition \ref{finiteprop}).
Suppose $f:X\ra S$ is
a finite surjective morphism between integral schemes,  and set 
$$
d=\deg(f) = [k(X):k(S)].
$$  
Then
$f_*$ is induced by pushforward of cycles, which in this
case is the map sending $[\{g\} \times X]$ to $d[\{g\} \times S]$.
Thus $f_*: \hig(X,G;\zz)\ra \hig(S,G;\zz)$ is simply multiplication by $d$.

\begin{prop} \label{algclosedprop}
Let $k$ be an algebraically closed field and let $A$ be an abelian
group.
Then we have natural (in $G$ and $A$) isomorphisms
\begin{align*}
\hig(\Spec(k),G;A) & \cong H_i(G(k), A) \\
\hgi(\Spec(k),G;A) & \cong H^i(G(k), A), 
\end{align*}
for all $i\ge 0$, 
where $G(k)$ is the discrete group of $k$-rational points of $G$.
\end{prop}

\begin{proof} This follows from the observation that
$C_0(G^i/ \Spec(k))$ is the free abelian group on 
$k$-points of $G^i$.
\end{proof}

If $k$ is not algebraically closed, then the situation is more complicated.

\begin{prop}\label{field} Let $k$ be a field and $\kbar$ be an 
  algebraic closure of $k$.  Let $\Gamma$ denote the absolute Galois
group $\text{\em Gal}(\kbar/k)$.  Then for all $i\ge 0$,
$$
\hig(\Spec(k),G;\zz)
\cong
H_i(BG(\kbar)/\Gamma).
$$
\end{prop}

\begin{proof} Observe that $C_0(G^i/\Spec(k))$ is the
  free abelian group on the closed points of $G^i$, which are 
in one-to-one correspondence
with the orbits of the $\Gamma$ action on $G^i(\kbar)$. We thus have
$$
C_0(G^i/\Spec(k)) \cong
\zz\{G^i(\kbar)/\Gamma\}.
$$
\end{proof}

This point of view allows us to make some calculations in Section \ref{invariant}. 

Observe that a map $X\ra G^i$ naturally defines an element of
$C_0(G^i\times X/X)$, and hence there is a map of chain complexes
$$
\zz \Hom(X, BG) \to
C_0(BG \times X/ X). 
$$
For example, if $X = \Spec(R)$ is affine, then
a morphism $X\ra G^i$ is simply an
$R$-point of $G^i$ and hence we have the map of chain complexes
$$
C_\bullet(G(R)) \lra C_0(BG\times X/X).
$$

\begin{prop} Let $X=\Spec(R)$ be an affine scheme.
Then  there is a  natural  map
$$
H_i(G(R),\zz)\lra \hig(X,G;\zz)
$$ 
for each $i\ge 0$. \hfill
$\qed$
\end{prop}

These maps are isomorphisms in the contexts of Propositions
\ref{finiteprop} and \ref{algclosedprop}, but in general 
they need not be 
either injective 
or surjective. In fact, we have the following result.

\begin{prop}\label{h1gmspeck}  There is an isomorphism
$$
\cH_1(\text{\em Spec}(\zr),\gm;\zz) \cong \zr^\times_{>0}
$$
and the map $H_1(\gm(\zr);\zz)\ra \cH_1(\text{\em Spec}(\zr),\gm;\zz)$
is surjective with kernel $\{\pm 1\}$.
\end{prop}

\begin{proof} 
The abelian group $\cH_1(\text{Spec}(\zr),\gm;\zz)$ is generated by
classes of non-zero complex numbers, $[z]$ for $z \in \C^\times$, 
modulo the relations $[z] = [\overline{z}]$
and $[zw] = [z] + [w]$, for all $z,w \in \C^\times$.
The homomorphism  $\zz(\C^\times) \to \zr^\times_{>0}$ induced by sending $z \in
\C^\times$ 
to $\|z\|$ annihilates these relations and hence induces a 
map
$$
\cH_1(\text{Spec}(\zr),\gm;\zz) 
\to
\zr^\times_{>0}.
$$
Likewise, these relations show that the
function $\zr^\times_{>0} \to \cH_1(\text{Spec}(\zr),\gm;\zz)$ given by
$r \mapsto [r]$ is a homomorphism. The composition 
$$\zr^\times_{>0} \to \cH_1(\text{Spec}(\zr),\gm;\zz) \to
\zr^\times_{>0}$$ is clearly the identity.

For $w \in \C^\times$, we have $w = z^2$, for some $z$, and
hence $[w] = [z^2] = 2[z] = [z] + [\overline{z}] = [z\overline{z}] =
[\|w\|]$. This shows that the composition 
$$\cH_1(\text{Spec}(\zr),\gm;\zz) \to
\zr^\times_{>0} \to 
\cH_1(\text{Spec}(\zr),\gm;\zz)$$ is also the identity.

Finally, the composition $\R^\times \cong H_1(\gm(\zr); \zz) \to
\cH_1(\text{Spec}(\zr),\gm;\zz) \to \R^\times_{>0}$ is given by $r
\mapsto |r|$.
\end{proof}

\section{The spectral sequence and Friedlander's Isomorphism Conjecture}\label{ss}
Let $\Delta^q$ denote the linear hypersurface in ${\mathbb A}^{q+1}$
defined by the equation $t_0+t_1+\cdots +t_q=1$.  There are obvious
coface and codegeneracy maps
\begin{eqnarray*}
\delta^i & : & \Delta^{q-1}\lra \Delta^q \\
\sigma^i & : & \Delta^q \lra \Delta^{q-1}
\end{eqnarray*}
making $\Delta^\bullet$ a cosimplicial scheme.

Fix $s\ge 0$.  Applying the functor $C_0(G^s\times -/-)$ to
$\Delta^\bullet$ yields a simplicial abelian group
$C_0(G^s\times\Delta^\bullet/\Delta^\bullet)$.  We therefore have a
double complex $A_{\bullet\bullet}$ with
$$
A_{s,t} = C_0(G^s\times\Delta^t/\Delta^t).
$$ 
Applying the functor $\text{Hom}(-,\zn)$ to $A_{\bullet\bullet}$
yields a double cochain complex $E^{\bullet\bullet}$ with
$$
E^{s,t} = \text{Hom}(C_0(G^s\times\Delta^t/\Delta^t),\zn).
$$

As usual, we have two first quadrant spectral sequences converging to
the cohomology of the total complex.  Taking horizontal cohomology
first, we obtain
$$
E_1^{s,t} = {\mathcal H}^s(\Delta^t,G;\zn).
$$

Using the work of Suslin-Voevodsky and the second spectral sequence,
we may identify the abutment.  If
${\mathcal F}$ is a presheaf with transfers, define a presheaf
${\mathcal F}_t$ by ${\mathcal F}_t(X) = {\mathcal
F}(X\times\Delta^t)$.  Then ${\mathcal F}_\bullet$ is a simplicial
presheaf of abelian groups. Also, we let $C_\bullet({\mathcal F})$
denote the simplicial abelian group 
${\mathcal F}(\Delta^\bullet)$.
There are obvious maps of simplicial presheaves $C_\bullet({\mathcal
F})\ra {\mathcal F}_\bullet$ and ${\mathcal F}\ra {\mathcal
F}_\bullet$ (where $C_\bullet(\cF)$ is regarded as a degree-wise
constant presheaf).  Corollary 7.7 of \cite{susvoe} asserts that both
of these maps induce isomorphisms upon applying the functor
$\text{Ext}^\bullet_{\text{qfh}}(-,\zn)$, provided $k$ is
algebraically closed and 
$n$ is relatively prime to the exponential characteristic of $k$.
Moreover, these qfh Ext groups coincide with the \'etale Ext groups
$\text{Ext}^\bullet_{\text{\'et}}(-,\zn)$ by \cite[10.10]{susvoe}.

Now consider the case ${\mathcal F}=C_0(G^s\times -/-)$.  Taking
vertical homology in the double complex yields $s$-th column
$$
\text{Ext}^\bullet_{\text{qfh}}(C_0(G^s\times -/-),\zn) =
\text{Ext}^\bullet_{\text{\'et}}(C_0(G^s\times -/-),\zn).
$$ In turn,
the latter is
$\text{Ext}^\bullet_{\text{\'et}}(\zz\text{Hom}(-,G^s),\zn)$ since
$C_0(G^s\times -/-)[1/p]$ is the qfh sheafification of
$\zz[1/p]\text{Hom}(-,G^s)$.  But these groups are precisely the
\'etale cohomology groups $H^\bullet_{\text{\'et}}(G^s,\zn)$
(Corollary $7.8$ of \cite{susvoe}).  What we conclude, then, is that
the second spectral sequence is simply the usual spectral sequence for
computing the \'etale cohomology of the simplicial scheme $BG_k$:
$$E_1^{s,t} = H^t_{\text{\'et}}(G^s,\zn) \Longrightarrow H^{s+t}(BG_k,\zn)$$
(see \cite{friedbook}, p.~16).  We have thus proved the following result.

\begin{theorem} \label{ssthm}  
Let $k$ be an algebraically closed field and $n$ a positive integer relatively
prime to the exponential characteristic of $k$.
Then there is a first quadrant spectral sequence
\begin{equation} \label{M3}
E_1^{s,t} = {\mathcal H}^s(\Delta^t,G;\zn) \Longrightarrow
H_{\text{\em \'et}}^{s+t}(BG_k,\zn).
\end{equation}
\hfill $\qed$
\end{theorem}

\begin{remark}
The results of \cite{susvoe} are stated only
for fields of characteristic zero.  Work of de Jong \cite{dejong}
allows us to replace this assumption with 
the assumption that 
$n$ is relatively
prime to the exponential characteristic of $k$.
\end{remark}

Observe that we have isomorphisms 
$$
E_1^{s,0} = \cH^s(\Spec(k), G; \zn)
\cong H^s(G(k), \zn) = H^s(BG(k), \zn)
$$ 
by Proposition \ref{algclosedprop}, and under this isomorphism, the edge
homomorphism
$$
H_{\text{\'et}}^{s}(BG_k,\zn) \to E_1^{s,0} \cong H^s(BG(k), \zn)
$$ 
is the map conjectured to be an isomorphism by
Friedlander. Consequently, we have:

\begin{cor} Friedlander's Isomorphism Conjecture for an algebraically closed
field $k$, an algebraic group $G$ defined over $k$, and a positive integer $n$
relatively prime to the exponential characteristic of $k$ is
equivalent to the assertion that the map of complexes
\begin{equation} \label{M1}
\Tot\left(\Hom(C_0(BG \times \Deldot/\Deldot),\zn)\right) \to
\Hom(C_0(BG/ \Spec(k)),\zn)
\end{equation}
is a quasi-isomorphism.
\end{cor}

For example, one might optimistically conjecture that for each fixed
$t \geq 0$ the map
$$
\cH^i(\Delta^t, G; \zn) \to \cH^i(\Spec(k), G; \zn)
$$
is an isomorphism, a result which would imply 
that (\ref{M1}) is a quasi-isomorphism and hence the Isomorphism
Conjecture. Though we have no counter-example to this stronger
assertion, it seems unlikely to hold.

The evident generalization to arbitrary coefficients of
the above assertion --- i.e., the assertion 
that the map of complexes
\begin{equation} \label{M2}
C_0(BG/\Spec(k))\to
\Tot\left(C_0(BG \times \Deldot/ \Deldot )\right) 
\end{equation}
is a quasi-isomorphism
 --- turns out to be false: Take $G = \gm$ and for each fixed $t$ 
let $C_0(\gm^{\wedge \bullet} \times \Delta^t/\Delta^t)$ 
denote the normalized chain
complex obtained from 
$C_0(B\gm \times \Delta^t/ \Delta^t)$ by modding out
degeneracies.
Then the bicomplexes
$C_0(B \gm^{\wedge \bullet} \times \Deldot/\Deldot)$ 
and
$C_0(B \gm \times \Deldot/\Deldot)$ are quasi-isomorphic, and by
\cite{VSF}, we have that for each fixed $s$, 
the complex $C_0(\gm^{\wedge s} \times \Deldot/\Deldot)$ computes
the weight $s$ motivic cohomology of the field $k$:
$$
h_i \left(C_0(\gm^{\wedge s} \times \Deldot/\Deldot)\right) 
\cong 
H_\cM^{s-i}(k, \zz(s)).
$$
It follows that $h_2$ of the total complex associated to 
$C_0(\gm^\bullet \times \Deldot/\Deldot)$ 
is the second  $K$-group of the field $k$: 
$$
h_2(C_0(\gm^\bullet \times \Deldot/ \Deldot)) \cong
K_2(k) = \bigwedge\nolimits^2(k^\times)/\left< a \wedge (1-a) \, | \,
a, 1-a \in k^\times \right>.
$$
By contrast, we have $h_2(C_0(B\gm/\Spec(k))) \cong
\bigwedge^2(k^\times)$, so that (\ref{M2})
has a kernel at $h_2$ (given by the Steinberg relations).

The map (\ref{M1}) is a quasi-isomorphism for $G = \gm$ (with $\zn$
coefficients), however, since the Isomorphism Conjecture holds for
this algebraic group (see Corollary \ref{h4} below). 
This is connected to the fact that the motivic
cohomology of $k$ with $\zn$ coefficients coincides with the singular
cohomology of a point:
$$
h_i(C_0(\gm^{\wedge s} \times \Deldot/\Deldot))
\cong
H_\cM^{s-i}(k, \zn(s)) \cong 
\begin{cases} 
\zn, & \text{if $s=i$,} \\
0, & \text{otherwise.} 
\end{cases}
$$

The spectral sequence also allows us to deduce the following criterion
for verifying Friedlander's Isomorphism Conjecture:

\begin{theorem}\label{isoconjthm}  
Let $G$ be an algebraic group over an algebraically closed field $k$
and let $n$ be a positive integer relatively prime to the exponential
characteristic of $K$.
Assume that for all smooth $k$-varieties $X$ and closed points $x$ of $X$,
the map
\begin{equation}\label{henseliso}
H_i(G(k),\zn) \lra H_i(G(\Rbar),\zn)
\end{equation}
is an isomorphism
for all $i\le r$, where $\Rbar$ is the absolute integral closure of $R
= \cO_{X,x}^{hen}$, the henselization of the local ring of $X$ at $x$.
Then the map
$$
H^i_{\text{\em \'et}}(BG_k,\zn)\lra H^i(BG(k),\zn)
$$ is an isomorphism for all
$i\le r$ and, moreover, 
there is an exact sequence
$$
0\lra H^{r+1}_{\text{\em \'et}}(BG_k,\zn)\lra H^{r+1}(BG(k),\zn) \stackrel{d}{\lra}
{\mathcal H}^{r+1}(\Delta^1,G;\zn),
$$
where $d$ is the difference of the maps induced by the two face maps $\Delta^0\ra \Delta^1$.
\end{theorem}

\begin{proof} 
The first implication is well-known under the stronger hypothesis that the map
(\ref{henseliso}) is an isomorphism with $\Rbar$ replaced by $R$ (see e.g.~\cite{jardine}).
Our machinery allows us to use this weaker hypothesis.  Consider the spectral sequence
$$
E_1^{s,t} = {\mathcal H}^s(\Delta^t,G;\zn)\Longrightarrow
H^{s+t}_{\text{\'et}}(BG_k,\zn)
$$ 
of Theorem \ref{ssthm}. Since $\Spec(\Rbar)$ is a limit of qfh
neighborhoods of $x \in X$, our hypothesis implies that the functors
$\cH_i(-, G; \zn)$ are locally constant for the qfh topology.
If ${\mathcal F}$ is a
presheaf on the category of schemes over $k$, let $\cF_{\text{qfh}}$ 
denote the
sheafification of ${\mathcal F}$ in the qfh topology.
We then have the following chain of natural
isomorphisms (where Ab denotes the category of abelian groups):
\begin{eqnarray*}
\text{Ext}^\bullet_{\text{Ab}}(\hig(\Spec(k),G;\zn),\zn) & \cong &
   \text{Ext}^\bullet_{\text{qfh}}(\hig(-,G;\zn)_{\text{qfh}},\zn) \\
   & \cong &
   \text{Ext}^\bullet_{\text{qfh}}(\hig(-\times\Delta^\bullet,G;\zn)_\text{qfh},\zn)
   \\ & \cong &
   \text{Ext}^\bullet_{\text{qfh}}(\hig(\Delta^\bullet,G;\zn),\zn) \\
   & \cong &
   \text{Ext}^\bullet_{\text{Ab}}(\hig(\Delta^\bullet,G;\zn),\zn)
\end{eqnarray*}
where the first holds because the functor $\hig(-,G;\zn)$ is locally constant in the qfh
topology, the middle two follow
from Theorem 7.6 of \cite{susvoe}, and the last holds by definition.  Since the functors $\hig(-,G;\zn)$ are
$n$-torsion, the map $\hig(\Delta^\bullet,G;\zn)\ra \hig(\Spec(k),G;\zn)$ must be a weak
equivalence for $i\le r$ and hence the same is true of the map $$\hgi(\Spec(k),G;\zn)\lra \hgi(\Delta^\bullet,G;\zn)$$  
for $i\le r$.  Thus, for all $s\le r$,
$$E_2^{s,t} = \begin{cases}
                H^s(BG(k),\zn) & t=0 \\
                0 & t>0.
             \end{cases}$$
So we see that the map $H^i_{\text{\'et}}(BG_k,\zn)\ra H^i(BG(k),\zn)$
is an isomorphism for $i\le r$.  Moreover, we have $E_\infty^{r+1,0} =
E_2^{r+1,0}$ and therefore we obtain the exact sequence
$$
0\lra H^{r+1}_{\text{\'et}}(BG_k,\zn)\lra H^{r+1}(BG(k),\zn)
\stackrel{d}{\lra} {\mathcal H}^{r+1}(\Delta^1,G;\zn).
$$ 
\end{proof}

\begin{cor}\label{h4}  The natural map
$$
H^i_{\text{\em \'et}}(BG_k,\zn) \lra H^i(BG(k),\zn)$$ is an isomorphism in the 
following cases:
\begin{enumerate}
\item $G$ finite, solvable, or the normalizer of a maximal torus in a
  reductive group;
\item $G=\glm$ in cohomological degrees $i\le 3$.
\end{enumerate}
Moreover, there is an exact sequence
$$0\lra H^4_{\text{\em \'et}}(BGL_m,\zn)\lra H^4(BGL_m(k),\zn)\lra {\mathcal H}^4(\Delta^1,GL_m;\zn).$$
\end{cor}

\begin{proof}  Let $X$ be a smooth $k$-variety, $x \in X$ a closed
  point, and $\Rbar$ the absolute integral closure of $\cO_{X,x}^{hen}$.

If $G$ is finite, then clearly 
$$
H_i(G(k),\zn) \stackrel{\cong}{\lra} H_i(G(\Rbar),\zn)
$$
for all $i\ge 0$.  If $G$ is solvable, then $G$ has a descending
central series whose graded quotients are either ${\mathbb G}_m$ or
${\mathbb G}_a$.  Clearly,
$$
H_i({\mathbb G}_a(k),\zn) \lra H_i({\mathbb G}_a(\Rbar),\zn)
$$ is an isomorphism for $i\ge 0$.  An
easy application of Hensel's lemma shows that the same is true for ${\mathbb G}_m$.  By
iterated use of the Hochschild--Serre spectral sequence we see that the map
$$
H_i(G(k),\zn)\lra H_i(G(\Rbar),\zn)
$$ is an isomorphism for all $i\ge 0$.  If $T$ is
a maximal torus in a reductive group $S$, then there is a short exact sequence
$$
1\lra T\lra N_S(T)\lra W\lra 1
$$  where $W$ is finite.  Again, the Hochschild--Serre
spectral sequence shows that the map
$$H_i(N_S(T)(k),\zn)\lra H_i(N_S(T)(\Rbar),\zn)$$ is an isomorphism for $i\ge 0$.  (All
the preceding facts may be found in \cite{jardine}.)

Finally, if $G=\glm$, then for $n$ prime to the characteristic of $k$ the map
$$H_i(G(k),\zn)\lra H_i(G(\Rbar),\zn)$$ is an isomorphism for $i\le 3$ (\cite{knudson}, p.~146).
Therefore, the isomorphism conjecture holds in this case as well and we obtain the above mentioned
exact sequence.
\end{proof}

\section{Calculations}\label{invariant}

For an arbitrary field $k$, 
Proposition \ref{field} relates the groups $\cH_\bullet(\Spec(k),
G;\zz)$ 
with a construction using a quotient of the action of the
absolute Galois group $\Gamma = \text{Gal}(\bark/k)$:
$$
\hig(\Spec(k),G;\zz) 
\cong 
H_i(BG(\kbar)/\Gamma).
$$

Now consider the field $k=\zr$.  We have $\overline{\zr}=\zc$ and $\text{Gal}(\overline{k}/k)
=\zz/2$.  If $G$ is an algebraic group over $\zr$, then we have an isomorphism
$${\mathcal H}_\bullet(\text{Spec}(\zr),G;A) = H_\bullet(BG(\zc)/\Gamma;A).$$
We therefore have the following result.

\begin{theorem}\label{specr} Let $G$ be an algebraic group over $\zr$ and let $A=\zq$ or
$\zn$, where $n$ is odd.  Then there is a canonical isomorphism
$$\hig(\text{\em Spec}(\zr),G;A) \stackrel{\cong}{\lra} H_i(G(\zc),A)^{\zz/2}$$
for each $i\ge 0$. 
\end{theorem}

\begin{proof}  This is a standard fact in the homology of quotient spaces, proved using the transfer
map.  See \cite{borel}, p.~38.
\end{proof}

If $A={\mathbb Z}/2$, the calculation is more difficult and we are only able to handle two
cases.  First, consider the case $G=\gm$.  We are interested in calculating
$H_\bullet(B\gm(\zc)/\Gamma;\zz/2)$.  Since homology and cohomology are dual with field
coefficients, it suffices to compute $H^\bullet(B\gm(\zc)/\Gamma;\zz/2)$.  While it is
possible to do this directly using various standard techniques, we proceed as follows.

Suppose $Q$ is a finite group acting on a CW-complex $X$ in such a way that if an element
of $Q$ fixes a cell of $X$, then it fixes it pointwise.  There is associated to this a 
cohomology theory, called (ordinary) Bredon cohomology,
$H^\bullet_Q(X; M)$, where $M$ is a Mackey functor, with the property
that
$$
H^\bullet_Q(X;\underline{\zz}) \cong H^\bullet(X/Q;\zz).
$$
Here, $\underline{A}$ is the constant Mackey functor associated to $A$.  For a thorough 
discussion of Bredon cohomology, we refer the reader to \cite{may}.  What is relevant for us
is the following result.

\begin{prop}\label{fixediso} Suppose $f:X\ra Y$ is a map of $Q$-CW-complexes such that for
each subgroup $H$ of $Q$ the induced map $f^H:X^H\ra Y^H$ of fixed point spaces induces an
isomorphism
$$H_\bullet(X^H;\zz/p) \lra H_\bullet(Y^H;\zz/p).$$  Then the induced map
$$f^*:H^\bullet_Q(Y;\underline{\zz/p})\lra H^\bullet_Q(X;\underline{\zz/p})$$
is an isomorphism.
\end{prop}

\begin{proof} See \cite{may}, p.~26.
\end{proof}

This suggests the following.

\begin{guess}\label{equivariantisoconj}  Let $G$ be an algebraic group over $\zr$. Let $\Gamma=
\text{\em Gal}(\zc/\zr)$, and let $p$ be a prime number.  Then the identity map
$G(\zc)\ra G(\zc)^{\text{\em top}}$ induces an isomorphism
$$
H^\bullet_\Gamma(BG(\zc)^{\text{\em top}};\underline{\zz/p})\lra
H^\bullet_\Gamma(BG(\zc);\underline{\zz/p}).
$$ 
\end{guess}

\begin{prop}\label{isoconjimpliesequiv} Friedlander's isomorphism conjecture for $G(\zc)$ and $G(\zr)$
implies the equivariant isomorphism conjecture for $G(\zc)$.
\end{prop}

\begin{proof}  Note that for $\Gamma=\zz/2$, the only subgroups are
  $\Gamma$ and the trivial subgroup 
$\{1\}$, and the corresponding fixed point spaces are $BG(\zr)$ and
$BG(\zc)$, respectively.  Thus, if the maps
$$H_\bullet(BG(\zr);\zz/p) \lra H_\bullet(BG(\zr)^{\text{top}};\zz/p)$$
and
$$H_\bullet(BG(\zc);\zz/p) \lra H_\bullet(BG(\zc)^{\text{top}};\zz/p)$$
are both isomorphisms, then Proposition \ref{fixediso} implies that the map
$$H^\bullet_\Gamma(BG(\zc)^{\text{top}};\underline{\zz/p})\lra
H^\bullet_\Gamma(BG(\zc);\underline{\zz/p})$$ 
is an isomorphism.
\end{proof}

\begin{cor}  Let $G$ be a solvable Lie group.  Then the equivariant
  isomorphism conjecture 
holds for $G$. 
\end{cor}

\begin{proof}  Friedlander's isomorphism conjecture holds for $G(\zc)$
  and $G(\zr)$ \cite{milnor}. 
\end{proof}

In the case $G=\gm$, we see that there is an isomorphism
$$H^\bullet_\Gamma(B\gm(\zc)^{\text{top}};\underline{\zz/2}) \lra
H^\bullet_\Gamma(B\gm(\zc);\underline{\zz/2}).$$ We are trying to
calculate the latter; these are the groups
$\cH^\bullet(\text{Spec}(\zr),\gm;\zz/2)$.  Note, however, that
$B\gm(\zc)^{\text{top}}$ is (equivariantly) homotopy equivalent to
$\zc \bP^\infty$.  So we must compute the groups
$$H^\bullet_\Gamma(\zc \bP^\infty;\underline{\zz/2}).
$$

Associated to a $\zz/2$-CW-complex $X$ is its $\zz/2$-equivariant
cohomology, which forms a bigraded ring
$$
H^{\bullet,\bullet}(X;\underline{A}),
$$ 
and extends the Bredon cohomology ring in the sense that 
we have
$$
H^{k,0}(X;\underline{\zz}) \cong H^k_{\zz/2}(X;\underline{\zz}).
$$
We shall not need the detailed definition of this theory, but the
interested reader may consult \cite{may}.  What is important for us is
the following result.

\begin{prop}\label{cpinftycoho} The cohomology of $\zc \bP^\infty$ is given by
$$H^{\bullet,\bullet}(\zc \bP^\infty;\underline{\zz/2}) \cong
  H^{\bullet,\bullet}(pt;\underline{\zz/2})[c],
$$ 
where $\deg c = (2,1)$.
\end{prop}

\begin{proof} See \cite{dugger}, 5.4, p.~18.
\end{proof}

The $\zz/2$-equivariant cohomology of a point with $\underline{\zz}$
coefficients has been calculated
(see \cite{dugger}, Appendix B), and from it one deduces
$$
H^{p,q}(pt;\underline{\zz/2}) \cong \begin{cases}
\zz/2 & \text{$q\ge p\ge 0$ or $q+2\le p\le 0$} \\
                                           0 & \text{otherwise}.
                                           \end{cases}
$$
Moreover, 
there is a commutative product on $H^{\bullet,\bullet}(pt)$ with the
property that the product of any element in degree $(p,q)$, $q\ge p\ge
0$ with an element of degree $(i,j)$, $j+2\le i\le 0$ is zero.

\begin{prop}\label{cpquotient} The cohomology groups $H^{s,0}(\zc
  \bP^\infty;\underline{\zz/2})$ satisfy $H^{1,0} = 0$ and for $k\ge
1$, $H^{2k,0} \cong H^{2k+1,0} \cong (\zz/2)^{k-1}$.
\end{prop}

\begin{proof}  Observe that elements of $H^{s,0}(\zc \bP^\infty)$
  arise only as products of powers of the 
generator $c$ (of degree $(2,1)$) with elements of
$H^{\bullet,\bullet}(pt)$ in degrees $(p,q)$ with $q+2\le p\le 0$.  As
there are no elements in degrees $(-1,-1)$, $(0,-1)$, or $(-1,-2)$, we
see that $H^{s,0}$ vanishes for $s=1,2,3$.  Denote the generator of
$H^{s,t}(pt)$ by $x_{(s,t)}$.  Then one sees easily that for $k\ge 2$,
we have
$$\begin{array}{c|c}
  \text{group} & \text{generators}  \\ \hline
                &                        \\
  H^{2k,0} & x_{(0,-k)}c^k, x_{(-2,-(k+1))}c^{k+1},\dots ,x_{(4-2k,2-2k)}c^{2k-2}  \\
           &          \\
  H^{2k+1,0} & x_{(-1,-(k+1))}c^{k+1},x_{(-3,-(k+2))}c^{k+2},\dots ,x_{(3-2k,1-2k)}c^{2k-1}
  \end{array}$$
This completes the proof.
\end{proof}

Recall that we have isomorphisms
$$
H^\bullet_\Gamma(\zc \bP^\infty;\underline{\zz/2})\cong H^\bullet(B\gm(\zc)/\Gamma;\zz/2)
\cong \cH^\bullet(\text{Spec}(\zr),\gm;\zz/2).
$$ 
Proposition \ref{cpquotient} therefore
gives us the following result.

\begin{theorem}\label{specrgm}  For all $k\ge 1$,
$$\cH^{2k}(\text{\em Spec}(\zr),\gm;\zz/2)\cong \cH^{2k+1}(\text{\em Spec}(\zr),\gm;\zz/2)\cong (\zz/2)^{k-1}.$$
The same is therefore true for $\cH_{2k}$ and $\cH_{2k+1}$.  \hfill $\qed$
\end{theorem}

Now suppose $G$ is a unipotent group over $\zr$.  Let $X=BG(\zc)$ and $Y=BG(\zr)=X^\Gamma$.  According to \cite{may},
p.~35, there is a long exact sequence (with $\zz/2$ coefficients)
$$\tilde{H}^n((X/Y)/\Gamma) \ra H^n(X)\ra \tilde{H}^n((X/Y)/\Gamma)\oplus H^n(Y)\ra \tilde{H}^{n+1}((X/Y)/\Gamma).$$
Note, however, that
$$(X/Y)/\Gamma = (X/\Gamma)/Y$$
and 
$$\tilde{H}^\bullet((X/\Gamma)/Y)\cong H^\bullet(X/\Gamma,Y).$$
Thus, this sequence becomes
$$H^n(X/\Gamma,Y)\ra H^n(X) \ra H^n(X/\Gamma,Y)\oplus H^n(Y)\ra H^{n+1}(X/\Gamma,Y).$$
Since $\tilde{H}^n(X) = \tilde{H}^n(Y) = 0$ in this case, we see that $H^n(X/\Gamma,Y)=0$
for all $n\ge 0$.  The long exact sequence of the pair $(X/\Gamma,Y)$ then shows that
$$\tilde{H}^n(X/\Gamma)=0$$
for all $n\ge 0$.  We therefore have the following result.

\begin{theorem} Let $G$ be a unipotent group over $\zr$.  Then for all $i>0$, $\hgi(\text{\em Spec}(\zr),G;\zz/2) = 0$.
The same is therefore true for $\hig$. \hfill $\qed$
\end{theorem}

\end{document}